# Fixed point theorems of soft contractive mappings

**Murat I. Yazar, Cigdem Gunduz (Aras), Sadi Bayramov**

Department of Mathematics, Kafkas University, Kars, 36100-Turkey
Department of Mathematics, Kocaeli University, Kocaeli, 41380-Turkey
Department of Mathematics, Kafkas University, Kars, 36100-Turkey

Email: miy248@yahoo.com, carasgunduz@gmail.com, baysadi@gmail.com

**Abstract:** The first aim of this paper is to examine some important properties of soft metric spaces. Second is to introduce soft continuous mappings and investigate properties of soft continuous mappings. Third is to prove some fixed point theorems of soft contractive mappings on soft metric spaces.

**Keywords:** soft metric space, soft continuous mapping, soft contractive mapping, fixed point theorem.

## 1. Introduction

In the year 1999, Molodtsov [11] initiated a novel concept of soft sets theory as a new mathematical tool for dealing with uncertainties. A soft set is a collection of approximate descriptions of an object. Soft systems provide a very general framework with the involvement of parameters. Since soft set theory has a rich potential, applications of soft set theory in other disciplines and real life problems are progressing rapidly.

Maji et al. [8,9] worked on soft set theory and presented an application of soft sets in decision making problems. Chen [2] introduced a new definition of soft set parameterization reduction and a comparison of it with attribute reduction in rough set theory. Many researchers contributed towards many structure on soft set theory. [1,2]

M.Shabir and M.Naz [12] presented soft topological spaces and they investigated some properties of soft topological spaces. Later, many researches about soft topological spaces were studied in [6,7,10,12,13]. In these studies, the concept of soft point is expressed by different approaches. In the study we use the concept of soft point which was given in [4,13].

It is known that there are many generalizations of metric spaces: Menger spaces, fuzzy metric spaces, generalized metric spaces, abstract (cone) metric spaces, or K-metric and K-normed spaces etc. Recently in [4,5] S.Das and S.K.Samanta introduced a different notion of soft metric space by using a different concept of soft point and investigated some basic properties of these spaces.

A number of authors have defined contractive type mapping on a complete metric space $X$ which are generalizations of the well-known Banach contraction, and which have the property that each such mapping has a unique fixed point [14,15]. The fixed point can always be found by using Picard iteration, beginning with some initial choice $x_0 \in X$.

In the present paper, we first give, as a preliminaries, some well-known results in soft set theory. Firstly, we examine some important properties of soft metric spaces defined in [4]. Secondly, we investigate properties of soft continuous mappings on soft metric spaces. Finally, we introduced soft contractive mappings on soft metric spaces and prove some fixed point theorems of soft contractive mappings.

## 2. Preliminaries

**Definition 2.1. [11]** Let $X$ be an initial universe set and $E$ be a set of parameters. A pair $(F,E)$ is called a soft set over $X$ if only if $F$ is a mapping from $E$ into the set of all subsets of the set $X$, i.e., $F: E \to P(X)$, where $P(X)$ is the power set of $X$.

**Definition 2.2. [9]** The intersection of two soft sets $(F,A)$ and $(G,B)$ over $X$ is the soft set $(H,C)$, where $C = A \cap B$ and $\forall \varepsilon \in C$, $H(\varepsilon) = F(\varepsilon) \cap G(\varepsilon)$. This is denoted by $(F,A) \tilde{\cap} (G,B) = (H,C)$.

**Definition 2.3. [9]** The union of two soft sets $(F,A)$ and $(G,B)$ over $X$ is the soft set, where $C = A \cup B$ and $\forall \varepsilon \in C$,
$$H(\varepsilon) = \begin{cases} F(\varepsilon), & \text{if } \varepsilon \in A - B \\ G(\varepsilon), & \text{if } \varepsilon \in B - A \\ F(\varepsilon) \cup G(\varepsilon), & \varepsilon \in A \cap B \end{cases}.$$
This relationship is denoted by $(F,A) \tilde{\cup} (G,B) = (H,C)$.

**Definition 2.4. [9]** A soft set $(F,A)$ over $X$ is said to be a null soft set denoted by $\Phi$ if for all $\varepsilon \in A$, $F(\varepsilon) = \emptyset$ (null set).

**Definition 2.5. [9]** A soft set $(F,A)$ over $X$ is said to be an absolute soft set, if for all $\varepsilon \in A, F(\varepsilon) = X$.

**Definition 2.6. [10]** The difference $(H,E)$ of two soft sets $(F,E)$ and $(G,E)$ over $X$, denoted by $(F,E) \setminus (G,E)$, is defined as $H(e) = F(e) \setminus G(e)$ for all $e \in E$.

**Definition 2.7. [10]** The complement of a soft set $(F,A)$ is denoted by $(F,A)^c$ and is defined by $(F,A)^c = (F^c, A)$, where $F^c: A \to P(X)$ is a mapping given by $F^c(\alpha) = X - F(\alpha), \forall \alpha \in A$.

**Definition 2.8. [3]** Let $\mathbb{R}$ be the set of real numbers and $B(\mathbb{R})$ be the collection of all non-empty bounded subsets of $\mathbb{R}$ and $E$ taken as a set of parameters. Then a mapping $F: E \to B(\mathbb{R})$ is called a soft real set. It is denoted by $(F,E)$. If specifically $(F,E)$ is a singleton soft set, then identifying $(F,E)$ with the corresponding soft element, it will be called a soft real number and denoted $\tilde{r}, \tilde{s}, \tilde{t}$ etc.
$\overline{0}, \overline{1}$ are the soft real numbers where $\overline{0}(e) = 0, \overline{1}(e) = 1$ for all $e \in E$, respectively.

**Definition 2.9. [3]** For two soft real numbers
(i) $\tilde{r} \leq \tilde{s}$ if $\tilde{r}(e) \leq \tilde{s}(e)$, for all $e \in E$;
(ii) $\tilde{r} \geq \tilde{s}$ if $\tilde{r}(e) \geq \tilde{s}(e)$, for all $e \in E$;
(iii) $\tilde{r} < \tilde{s}$ if $\tilde{r}(e) < \tilde{s}(e)$, for all $e \in E$;

(iv) $\tilde{r} > \tilde{s}$ if $\tilde{r}(e) > \tilde{s}(e)$, for all $e \in E$.

**Definition 2.10. [4,13]** A soft set $(P,E)$ over $X$ is said to be a soft point if there is exactly one $e \in E$, such that $P(e) = \{x\}$ for some $x \in X$ and $P(e') = \emptyset, \forall e' \in E/\{e\}$. It will be denoted by $\tilde{x}_e$.

**Definition 2.11. [4,13]** Two soft point $\tilde{x}_e, \tilde{y}_{e'}$ are said to be equal if $e = e'$ and $P(e) = P(e')$ i.e., $x = y$. Thus $\tilde{x}_e \neq \tilde{y}_{e'} \Leftrightarrow x \neq y$ or $e \neq e'$.

**Proposition 2.12. [13]** The union of any collection of soft points can be considered as a soft set and every soft set can be expressed as union of all soft points belonging to it; $(F, E) = \bigcup_{\tilde{x}_e \in (F,E)} \tilde{x}_e$.

**Definition 2.13. [12]** Let $\tau$ be the collection of soft sets over $X$, then $\tau$ is said to be a soft topology on $X$ if
(1) $\Phi, X$ belong to $\tau$
(2) the union of any number of soft sets in $\tau$ belongs to $\tau$
(3) the intersection of any two soft sets in $\tau$ belongs to $\tau$.

The triplet $(X, \tau, E)$ is called a soft topological space over $X$.

**Definition 2.14. [7]** Let $(X, \tau, E)$ be a soft topological space over $X$. Then soft interior of $(F, E)$, denoted by $(F, E)^\circ$, is defined as the union of all soft open sets contained in $(F, E)$.

**Definition 2.15. [7]** Let $(X, \tau, E)$ be a soft topological space over $X$. Then soft closure of $(F, E)$, denoted by $\overline{(F, E)}$, is defined as the intersection of all soft closed super sets of $(F, E)$.

**Definition 2.16. [7]** Let $(X, \tau, E)$ be a soft topological space over $X$. Then the soft boundary of soft set $(F, E)$ over $X$ is denoted by $\partial(F, E)$ and is defined as $\partial(F, E) = \overline{(F, E)} \tilde{\cap} \overline{(F, E)^c}$.

**Definition 2.17. [6]** Let $(X, \tau, E)$ and $(Y, \tau', E)$ be two soft topological spaces, $f : (X, \tau, E) \to (Y, \tau', E)$ be a mapping. For each soft neighborhood $(H, E)$ of $(f(x)_e, E)$, if there exists a soft neighborhood $(F, E)$ of $(x_e, E)$ such that $f((F, E)) \subset (H, E)$, then $f$ is said to be soft continuous mapping at $(x_e, E)$.
If $f$ is soft continuous mapping for all $(x_e, E)$, then $f$ is called soft continuous mapping.

Let $\tilde{X}$ be the absolute soft set i.e., $F(e) = X, \forall e \in E$, where $(F, E) = \tilde{X}$ and $SP(\tilde{X})$ be the collection of all soft points of $\tilde{X}$ and $\mathbb{R}(E)^*$ denote the set of all non-negative soft real numbers.

**Definition 2.18. [4]** A mapping $\tilde{d}: SP(\tilde{X}) \times SP(\tilde{X}) \to \mathbb{R}(E)^*$, is said to be a soft metric on the soft set $\tilde{X}$ if $d$ satisfies the following conditions:

(M1) $\tilde{d}(\tilde{x}_{e_1}, \tilde{y}_{e_2}) \tilde{\geq} \overline{0}$ for all $\tilde{x}_{e_1}, \tilde{y}_{e_2} \tilde{\in} \tilde{X}$,

(M2) $\tilde{d}(\tilde{x}_{e_1}, \tilde{y}_{e_2}) = \overline{0}$ if and only if $\tilde{x}_{e_1} = \tilde{y}_{e_2}$,

(M3) $\tilde{d}(\tilde{x}_{e_1}, \tilde{y}_{e_2}) = \tilde{d}(\tilde{y}_{e_2}, \tilde{x}_{e_1})$ for all $\tilde{x}_{e_1}, \tilde{y}_{e_2} \tilde{\in} \tilde{X}$,

(M4) For all $\tilde{x}_{e_1}, \tilde{y}_{e_2}, \tilde{z}_{e_3} \tilde{\in} \tilde{X}$, $\tilde{d}(\tilde{x}_{e_1}, \tilde{z}_{e_3}) \tilde{\leq} \tilde{d}(\tilde{x}_{e_1}, \tilde{y}_{e_2}) + \tilde{d}(\tilde{y}_{e_2}, \tilde{z}_{e_3})$.

The soft set $\tilde{X}$ with a soft metric $\tilde{d}$ on $\tilde{X}$ is called a soft metric space and denoted by $(\tilde{X}, \tilde{d}, E)$.

**Definition 2.19. [4]** Let $(\tilde{X}, \tilde{d}, E)$ be a soft metric space and $\tilde{\varepsilon}$ be a non-negative soft real number. $B(\tilde{x}_e, \tilde{\varepsilon}) = \{\tilde{y}_{e'} \tilde{\in} \tilde{X} : \tilde{d}(\tilde{x}_e, \tilde{y}_{e'}) \tilde{<} \tilde{\varepsilon}\} \subset SP(\tilde{X})$ is called the soft open ball with center $\tilde{x}_e$ and radius $\tilde{\varepsilon}$ and $B[\tilde{x}_e, \tilde{\varepsilon}] = \{\tilde{x}_e \tilde{\in} \tilde{X}; \tilde{d}(\tilde{x}_e, \tilde{y}_{e'}) \tilde{\leq} \tilde{\varepsilon}\} \subset SP(\tilde{X})$ is called the soft closed ball with center $\tilde{x}_e$ and radius $\tilde{\varepsilon}$.

**Definition 2.20. [4]** Let $(\tilde{X}, \tilde{d}, E)$ be a soft metric space and $(F, E)$ be a non-null soft subset of $\tilde{X}$ in $(\tilde{X}, \tilde{d}, E)$. Then $(F, E)$ is said to be soft open in $\tilde{X}$ with respect to $\tilde{d}$ if and only if all soft points of $(F, E)$ is interior points of $(F, E)$.

**Definition 2.21. [4]** Let $\{\tilde{x}_{\lambda, n}\}_n$ be a sequence of soft points in a soft metric space $(\tilde{X}, \tilde{d}, E)$. The sequence $\{\tilde{x}_{\lambda, n}\}_n$ is said to be convergent in $(\tilde{X}, \tilde{d}, E)$ if there is a soft point $\tilde{y}_\mu \tilde{\in} \tilde{X}$ such that $d(\tilde{x}_{\lambda, n}, \tilde{y}_\mu) \to \overline{0}$ as $n \to \infty$.

This means for every $\tilde{\varepsilon} \tilde{>} \overline{0}$, chosen arbitrarily, $\exists$ a natural number $N = N(\tilde{\varepsilon})$, such that $\overline{0} \tilde{\leq} d(\tilde{x}_{\lambda, n}, \tilde{y}_\mu) \tilde{<} \tilde{\varepsilon}$, whenever $n > N$.

**Theorem 2.22. [4]** Limit of a sequence in a soft metric space, if exist is unique.

**Definition 2.23. [4]** (Cauchy Sequence). A sequence $\{\tilde{x}_{\lambda, n}\}_n$ of soft points in $(\tilde{X}, \tilde{d}, E)$ is considered as a Cauchy sequence in $\tilde{X}$ if corresponding to every $\tilde{\varepsilon} \tilde{>} \overline{0}, \exists m \in N$ such that $d(\tilde{x}_{\lambda, i}, \tilde{x}_{\lambda, j}) \tilde{\leq} \tilde{\varepsilon}, \forall i, j \geq m$, i.e., $d(\tilde{x}_{\lambda, i}, \tilde{x}_{\lambda, j}) \to \overline{0}$ as $i, j \to \infty$.

**Definition 2.24. [4]** (Complete Metric Space). A soft metric space $(\tilde{X}, \tilde{d}, E)$ is called complete if every Cauchy Sequence in $\tilde{X}$ converges to some point of $\tilde{X}$. The soft metric space $(\tilde{X}, \tilde{d}, E)$ is called incomplete if it is not complete.

## 3. Soft Topology Generated by Soft Metric

In this section, we study some important results of soft metric spaces.

Let $\tilde{X}$ be the absolute soft set over $E$ and $\tilde{X}_\lambda$ be a family of soft points i.e., $\tilde{X}_\lambda = \{\tilde{x}_\lambda : x \in X\}$, for $\forall \lambda \in E$. Then there exists a bijective mapping between the soft set $\tilde{X}_\lambda$ and the set $X$.

If $\lambda \neq \mu \in E$, then $\tilde{X}_\lambda \cap \tilde{X}_\mu = \Phi$ and $\tilde{X} = \bigcup_{\lambda \in E} \tilde{X}_\lambda$.

Let $(\tilde{X}, \tilde{d}, E)$ be a soft metric space. It is clear that $(\tilde{X}_\lambda, \tilde{d}_\lambda, \{\lambda\})$ is a soft metric space, for $\lambda \in E$. Then by using the soft metric $\tilde{d}_\lambda$, we define a metric on $X$ as $d_\lambda(x, y) = \tilde{d}_\lambda(\tilde{x}_\lambda, \tilde{y}_\lambda)$. Note that $\lambda \neq \mu \in E$, then $d_\lambda$ and $d_\mu$ on $X$ are generally different metrics.

**Proposition 3.1.** Every soft metric space is a family of parameterized metric spaces.

The converse of the Proposition 3.1 may not be true in general. This is shown by the following example.

**Example 3.2.** Let $E = \mathbb{R}$ be a parameter set and $(X, d)$ be a metric space. We define the function $\tilde{d} : SP(\tilde{X}) \times SP(\tilde{X}) \to \mathbb{R}(E)^*$ by, $\tilde{d}(\tilde{x}_\lambda, \tilde{y}_\mu) = d(x, y)^{1+|\lambda - \mu|}$ for all $\tilde{x}_\lambda, \tilde{y}_\mu \in SP(\tilde{X})$. Then for all $\lambda \in E$, $d_\lambda$ is a metric on $X$. If $\tilde{d}(\tilde{x}_\lambda, \tilde{y}_\mu) = \overline{0}$ then this does not always mean that $\tilde{x}_\lambda = \tilde{y}_\mu$, so $\tilde{d}$ is not a soft metric on $\tilde{X}$.

**Proposition 3.3.** Let $(\tilde{X}, \tilde{d}, E)$ be a soft metric space and $\tau_{\tilde{d}}$ be a soft topology generated by the soft metric $\tilde{d}$. Then for every $\lambda \in E$, the topology $(\tau_{\tilde{d}})_\lambda$ on $X$ is the topology $\tau_{d_\lambda}$ generated by the metric $d_\lambda$ on $X$.
**Proof.** The proof is obvious.

**Lemma 3.4.** Let $(\tilde{X}, \tilde{d}, E)$ be a soft metric space. Then the following expressions are true.
(i) $\tilde{x}_\lambda \tilde{\in} \overline{(F, E)} \Leftrightarrow \tilde{d}(\tilde{x}_\lambda, (F, E)) = \overline{0}$;
(ii) $\tilde{x}_\lambda \tilde{\in} (F, E)^\circ \Leftrightarrow \tilde{d}(\tilde{x}_\lambda, (F, E)^c) > \overline{0}$;
(iii) $\tilde{x}_\lambda \tilde{\in} \partial(F, E) \Leftrightarrow \tilde{d}(\tilde{x}_\lambda, (F, E)) = \tilde{d}(\tilde{x}_\lambda, (F, E)^c) = \overline{0}$;

Note that if $(F, E)$ is a soft closed set in the soft metric space $(\tilde{X}, \tilde{d}, E)$ and $\tilde{x}_\lambda \tilde{\notin} (F, E)$, then there exists a soft open ball $B(\tilde{x}_\lambda, \tilde{\varepsilon})$ such that $B(\tilde{x}_\lambda, \tilde{\varepsilon}) \tilde{\cap} (F, E) = \Phi$.

**Theorem 3.8.** Every soft metric space is a soft normal space.
**Proof.** Let $(F_1, E)$ and $(F_2, E)$ be two disjoint soft closed sets in the soft metric space $(\tilde{X}, \tilde{d}, E)$. For every soft points $\tilde{x}_\lambda \tilde{\in} (F_1, E)$ and $\tilde{y}_\mu \tilde{\in} (F_2, E)$, we choose soft open balls $B(\tilde{x}_\lambda, \tilde{\varepsilon}_{\tilde{x}_\lambda})$ and $B(\tilde{y}_\mu, \tilde{\delta}_{\tilde{y}_\mu})$ such that $B(\tilde{x}_\lambda, \tilde{\varepsilon}_{\tilde{x}_\lambda}) \tilde{\cap} (F_2, E) = \Phi$ and $B(\tilde{y}_\mu, \tilde{\delta}_{\tilde{y}_\mu}) \tilde{\cap} (F_1, E) = \Phi$.

Thus, we have $(F_1, E) \tilde{\subset} \tilde{\cup} B\left(\tilde{x}_\lambda, (\tilde{\varepsilon}/3)_{\tilde{x}_\lambda}\right) = (U, E)$ and $(F_2, E) \tilde{\subset} \tilde{\cup} B\left(\tilde{y}_\mu, (\tilde{\delta}/3)_{\tilde{y}_\mu}\right) = (V, E)$.

We want to show that $(U, E) \tilde{\cap} (V, E) = \Phi$.

Assume that $(U, E) \tilde{\cap} (V, E) \neq \Phi$. Then there exists a soft point $\tilde{z}_\nu$ such that $\tilde{z}_\nu \tilde{\in} (U, E) \tilde{\cap} (V, E)$. Therefore, there exist soft open balls $B\left(\tilde{x}_\lambda, (\tilde{\varepsilon}/3)_{\tilde{x}_\lambda}\right)$ and $B\left(\tilde{y}_\mu, (\tilde{\delta}/3)_{\tilde{y}_\mu}\right)$ such that $\tilde{z}_\nu \tilde{\in} B\left(\tilde{x}_\lambda, (\tilde{\varepsilon}/3)_{\tilde{x}_\lambda}\right)$ and $\tilde{z}_\nu \tilde{\in} B\left(\tilde{y}_\mu, (\tilde{\delta}/3)_{\tilde{y}_\mu}\right)$. Here, we have $\tilde{d}(\tilde{x}_\lambda, \tilde{z}_\nu) \tilde{<} (\tilde{\varepsilon}/3)_{\tilde{x}_\lambda}$ and $\tilde{d}(\tilde{y}_\mu, \tilde{z}_\nu) \tilde{<} (\tilde{\delta}/3)_{\tilde{y}_\mu}$. If we get $\max\left\{(\tilde{\varepsilon}/3)_{\tilde{x}_\lambda}, (\tilde{\delta}/3)_{\tilde{y}_\mu}\right\} = (\tilde{\varepsilon}/3)_{\tilde{x}_\lambda}$, then we have $\tilde{d}(\tilde{x}_\lambda, \tilde{y}_\mu) \tilde{\leq} \tilde{d}(\tilde{x}_\lambda, \tilde{z}_\nu) + \tilde{d}(\tilde{z}_\nu, \tilde{y}_\mu) \tilde{<} (\tilde{\varepsilon}/3)_{\tilde{x}_\lambda} + (\tilde{\delta}/3)_{\tilde{y}_\mu} \tilde{<} \tilde{\varepsilon}_{\tilde{x}_\lambda}$ and so $\tilde{y}_\mu \tilde{\in} B(\tilde{x}_\lambda, \tilde{\varepsilon}_{\tilde{x}_\lambda})$ and $\tilde{y}_\mu \tilde{\in} (F_2, E)$ which contradicts with our assumption. Therefore, $(U, E) \tilde{\cap} (V, E) = \Phi$.

## 4. Soft Contractive Mappings

In this section we shall prove some fixed point theorems of soft contractive mappings.

Let $(\tilde{X}, \tilde{d}, E)$ and $(\tilde{Y}, \tilde{\rho}, E')$ be two soft metric spaces. The mapping $(f, \varphi): (\tilde{X}, \tilde{d}, E) \to (\tilde{Y}, \tilde{\rho}, E')$ is a soft mapping, where $f: X \to Y$, $\varphi: E \to E'$ are two mappings.

**Proposition 4.1.** For each soft point $\tilde{x}_\lambda \tilde{\in} \tilde{X}$, $(f, \varphi)(\tilde{x}_\lambda)$ is a soft point in $\tilde{Y}$.

**Proof.** Let $\tilde{x}_\lambda \tilde{\in} \tilde{X}$ be a soft point. Then
$$(f, \varphi)(\tilde{x}_\lambda)(e') = \bigcup_{e \in \varphi^{-1}(e')} f(x_\lambda(e)) = \bigcup_{\varphi(\lambda) = e'} f(x_\lambda(e)) = (f(x))_{\varphi(\lambda)}.$$

**Definition 4.2.** Let $(\tilde{X}, \tilde{d}, E)$ and $(\tilde{Y}, \tilde{\rho}, E')$ be two metric spaces and $(f, \varphi): (\tilde{X}, \tilde{d}, E) \to (\tilde{Y}, \tilde{\rho}, E')$ be a soft mapping. The mapping $(f, \varphi): (\tilde{X}, \tilde{d}, E) \to (\tilde{Y}, \tilde{\rho}, E')$ is soft continuous at the soft point $\tilde{x}_\lambda \in SP(\tilde{X})$, if for every soft open ball $B((f, \varphi)(\tilde{x}_\lambda), \tilde{\varepsilon})$ of $(\tilde{Y}, \tilde{\rho}, E')$ there exists a soft open ball $B(\tilde{x}_\lambda, \tilde{\delta})$ of $(\tilde{X}, \tilde{d}, E)$ such that $f\left(B(\tilde{x}_\lambda, \tilde{\delta})\right) \tilde{\subseteq} B((f, \varphi)(\tilde{x}_\lambda), \tilde{\varepsilon})$.

If $(f, \varphi)$ is soft continuous at every soft point $\tilde{x}_\lambda$ of $(\tilde{X}, \tilde{d}, E)$, then it is said to be soft continuous on $(\tilde{X}, \tilde{d}, E)$.

It can be expressed using $\varepsilon - \delta$ as follows;

The mapping $(f,\varphi):(\tilde{X},\tilde{d},E) \to (\tilde{Y},\tilde{\rho},E')$ is said to be soft continuous at the soft point $\tilde{x}_\lambda \in SP(\tilde{X})$, if for every $\tilde{\varepsilon} > 0$ there exists a $\tilde{\delta} > 0$ such that $\tilde{d}(\tilde{x}_\lambda, \tilde{y}_\mu) \tilde{<} \tilde{\delta}$ implies that $\tilde{\rho}((f,\varphi)(\tilde{x}_\lambda),(f,\varphi)(\tilde{y}_\mu)) \tilde{<} \tilde{\varepsilon}$.

**Theorem 4.3.** Let $(f,\varphi):(\tilde{X},\tilde{d},E) \to (\tilde{Y},\tilde{\rho},E')$ be a soft mapping. Then the following conditions are equivalent:

(1) $(f,\varphi):(\tilde{X},\tilde{d},E) \to (\tilde{Y},\tilde{\rho},E')$ is a soft continuous mapping,

(2) For each soft open set $(G,E')$ over $Y$, $(f,\varphi)^{-1}((G,E'))$ is a soft open set over $X$,

(3) For each soft closed set $(H,E')$ over $Y$, $(f,\varphi)^{-1}((H,E'))$ is a soft closed set over $X$,

(4) For each soft set $(F,E)$ over $X$, $(f,\varphi)(\overline{(F,E)}) \subset \overline{((f,\varphi)(F,E))}$,

(5) For each soft set $(G,E')$ over $Y$, $\overline{((f,\varphi)^{-1}(G,E'))} \subset (f,\varphi)^{-1}(\overline{(G,E')})$,

(6) For each soft set $(G,E')$ over $Y$, $f^{-1}((G,E')^\circ) \subset (f^{-1}(G,E'))^\circ$.

**Proof.** (1)$\Rightarrow$(2) Let $(f,\varphi)$ be a soft continuous and $(G,E') \subset (\tilde{Y},\tilde{\rho},E')$ be a soft open set. Consider the soft set $(f,\varphi)^{-1}(G,E')$. If $(f,\varphi)^{-1}(G,E') = \Phi$, then the proof is completed. Let $(f,\varphi)^{-1}(G,E') \neq \Phi$. In this case there exists at least one soft point $\tilde{x}_\lambda \in (f,\varphi)^{-1}(G,E')$. Then we have $(f,\varphi)(\tilde{x}_\lambda) \in (G,E')$. Since $(G,E')$ is soft open, there exists soft open ball $B((f,\varphi)(\tilde{x}_\lambda),\tilde{\varepsilon})$ such that $B((f,\varphi)(\tilde{x}_\lambda),\tilde{\varepsilon}) \subset (G,E')$ holds. Also since $(f,\varphi)$ is soft continuous, there exists a soft open ball $B(\tilde{x}_\lambda,\tilde{\delta})$ such that $(f,\varphi)(B(\tilde{x}_\lambda,\tilde{\delta})) \subset B((f,\varphi)(\tilde{x}_\lambda),\tilde{\varepsilon})$. Thus,

$$B(\tilde{x}_\lambda,\tilde{\delta}) \subset (f,\varphi)^{-1}(f,\varphi)(B(\tilde{x}_\lambda,\tilde{\delta})) \subset (f,\varphi)^{-1}B((f,\varphi)(\tilde{x}_\lambda),\tilde{\varepsilon}) \subset (f,\varphi)^{-1}(G,E').$$

Consequently, $(f,\varphi)^{-1}(G,E')$ is a soft open set.

(2)$\Rightarrow$(3) Let $(H,E')$ be any soft closed set over $Y$. Then $(H,E')^c$ is a soft open set. From (2), we have $((f,\varphi)^{-1}((H,E')))^c$ is a soft open set over $X$. Thus $(f,\varphi)^{-1}((H,E'))$ is a soft closed set.

(3)$\Rightarrow$(4) Let $(F,E)$ be a soft set over $X$. Since
$$(F,E) \subset (f,\varphi)^{-1}((f,\varphi)(F,E)) \text{ and } (f,\varphi)(F,E) \subset \overline{((f,\varphi)(F,E))},$$
we have $(F,E) \subset (f,\varphi)^{-1}((f,\varphi)(F,E)) \subset (f,\varphi)^{-1}(\overline{(f,\varphi)(F,E)})$. By part (3), since $(f,\varphi)^{-1}(\overline{(f,\varphi)(F,E)})$ is a soft closed set over $X$, $\overline{(F,E)} \subset (f,\varphi)^{-1}(\overline{(f,\varphi)(F,E)})$. Thus $(f,\varphi)(\overline{(F,E)}) \subset (f,\varphi)((f,\varphi)^{-1}\overline{(f,\varphi)(F,E)}) \subset \overline{(f,\varphi)(F,E)}$ is obtained.

(4)$\Rightarrow$(5) Let $(G,E')$ be a soft set over $Y$ and $(f,\varphi)^{-1}(G,E') = (F,E)$. By part (4), we have

$(f,\varphi)\left(\overline{(F,E)}\right) = (f,\varphi)\left(\overline{(f,\varphi)^{-1}(G,E')}\right) \subset \overline{(f,\varphi)\left((f,\varphi)^{-1}(G,E')\right)} \subset \overline{(G,E')}$. Then

$\overline{(f,\varphi)^{-1}(G,E')} = \overline{(F,E)} \subset (f,\varphi)^{-1}\left((f,\varphi)\overline{(F,E)}\right) \subset (f,\varphi)^{-1}\left(\overline{(G,E')}\right)$.

**(5)$\Rightarrow$(6)** Let $(G,E)$ be a soft set over $Y$. Substituting $(G,E')^c$ for condition in (5). Then $\overline{(f,\varphi)^{-1}((G,E')^c)} \subset (f,\varphi)^{-1}\left(\overline{(G,E')^c}\right)$. Since $(G,E')^\circ = \left(\overline{(G,E')^c}\right)^c$, then we have

$(f,\varphi)^{-1}\left((G,E')^\circ\right) = (f,\varphi)^{-1}\left(\left(\overline{(G,E')^c}\right)^c\right) = \left((f,\varphi)^{-1}\left(\overline{(G,E')^c}\right)\right)^c \subset$

$\left(\overline{(f,\varphi)^{-1}((G,E')^c)}\right)^c = \left(\overline{\left((f,\varphi)^{-1}(G,E')\right)^c}\right)^c = \left((f,\varphi)^{-1}(G,E')\right)^\circ$.

**(6)$\Rightarrow$(1)** Let $(G,E')$ be a soft open set over $Y$. Then since

$\left((f,\varphi)^{-1}(G,E')\right)^\circ \subset (f,\varphi)^{-1}(G,E') = (f,\varphi)^{-1}\left((G,E')^\circ\right) \subset \left((f,\varphi)^{-1}(G,E')\right)^\circ$,

$\left((f,\varphi)^{-1}(G,E')\right)^\circ = (f,\varphi)^{-1}(G,E')$ is obtained. This implies that $f^{-1}(G,E')$ is a soft open set.

**Definition 4.4.** The soft mapping $(f,\varphi):(\tilde{X},\tilde{d},E) \rightarrow (\tilde{Y},\tilde{\rho},E')$ is said to be soft sequentially continuous at the soft point $\tilde{x}_\lambda \tilde{\in} SP(\tilde{X})$ iff for every sequence of soft points $\{\tilde{x}_{\lambda_n}\}_n$ converging to the soft point $\tilde{x}_\lambda$ in the metric space $(\tilde{X},\tilde{d},E)$, the sequence $(f,\varphi)\left(\{\tilde{x}_{\lambda_n}\}_n\right)$ in $(\tilde{Y},\tilde{\rho},E')$ converges to a soft point $(f,\varphi)(\tilde{x}_\lambda) \tilde{\in} SP(\tilde{Y})$.

**Theorem 4.5.** Soft continuity is equivalent to soft sequential continuity in soft metric spaces.
**Proof.** Let $(f,\varphi):(\tilde{X},\tilde{d},E) \rightarrow (\tilde{Y},\tilde{\rho},E')$ be a soft continuous mapping and $\{\tilde{x}_{\lambda_n}\}_n$ be any sequence of soft points converging to the soft point $\tilde{x}_\lambda$ in the soft metric space $(\tilde{X},\tilde{d},E)$. Let $B((f,\varphi)(\tilde{x}_\lambda),\tilde{\varepsilon})$ be a soft open ball in $(\tilde{Y},\tilde{\rho},E')$. By continuity of $(f,\varphi)$ choose a soft open ball $B(\tilde{x}_\lambda,\tilde{\delta})$ containing $\tilde{x}_\lambda$ in $X$ such that $f\left(B(\tilde{x}_\lambda,\tilde{\delta})\right) \tilde{\subseteq} B((f,\varphi)(\tilde{x}_\lambda),\tilde{\varepsilon})$. Since $\{\tilde{x}_{\lambda_n}\}_n$ converges $\tilde{x}_\lambda$ there exists $n_0 \in \mathbb{N}$ such that $\{\tilde{x}_{\lambda_n}\}_n \tilde{\in} B(\tilde{x}_\lambda,\tilde{\delta})$ for all $n \geq n_0$. Therefore we have that $f\left(\{\tilde{x}_{\lambda_n}\}_n\right) \tilde{\in} f\left(B(\tilde{x}_\lambda,\delta)\right) \tilde{\subseteq} B((f,\varphi)(\tilde{x}_\lambda),\tilde{\varepsilon})$ for all $n \geq n_0$, as required.
Conversely, assume for contradiction that $(f,\varphi):(\tilde{X},\tilde{d},E) \rightarrow (\tilde{Y},\tilde{\rho},E')$ is soft sequential continuous but not soft continuous. Since $(f,\varphi)$ is not soft continuous at $\tilde{x}_\lambda$, there exists $\tilde{\varepsilon}_0 \tilde{>} \bar{0}$ such that for all $\tilde{\delta} \tilde{>} \bar{0}$ there exists $\tilde{y}_\mu \tilde{\in} SP(\tilde{X})$ such that $\tilde{d}(\tilde{x}_\lambda,\tilde{y}_\mu) \tilde{<} \tilde{\delta}$ and $\tilde{\rho}\left((f,\varphi)(\tilde{x}_\lambda),(f,\varphi)(\tilde{y}_\mu)\right) \tilde{>} \tilde{\varepsilon}_0$. For $n \geq 1 (n \in \mathbb{N})$, define $\tilde{\delta}_n = \frac{1}{n}$. For $n \geq 1$ we may choose $(\tilde{y}_\mu)_n \tilde{\in} SP(\tilde{X})$ such that $\tilde{d}((\tilde{y}_\mu)_n, \tilde{x}_\lambda) \tilde{<} \tilde{\delta}_n$ and $\tilde{\rho}((f,\varphi)(\tilde{x}_\lambda),(f,\varphi)(\tilde{y}_\mu)) \tilde{>} \tilde{\varepsilon}_0$ Therefore, by definition the sequence $\{(\tilde{y}_\mu)\}_n$ $(n \geq 1)$ converges to $\tilde{x}_\lambda$. However, by definition the sequence $\{(f,\varphi)(\tilde{y}_\mu)\}_n$ $(n \geq 1)$ does not converge to $(f,\varphi)(\tilde{x}_\lambda)$. That is, $(f,\varphi)$ is not soft sequentially continuous at $\tilde{x}_\lambda$.

**Definition 4.6.** Let $(\tilde{X}, \tilde{d}, E)$ be a soft metric space. A function $(f,\varphi):(\tilde{X},\tilde{d},E) \to (\tilde{X},\tilde{d},E)$ is called a soft contraction mapping if there exists a soft real number $\alpha \in \mathbb{R}, \bar{0} \leq \alpha < \bar{1}$ such that for every soft points $\tilde{x}_\lambda, \tilde{y}_\mu \in SP(X)$ we have $\tilde{d}((f,\varphi)(\tilde{x}_\lambda),(f,\varphi)(\tilde{y}_\mu)) \tilde{\leq} \alpha.\tilde{d}(\tilde{x}_\lambda,\tilde{y}_\mu)$.

**Proposition 4.7.** Every soft contraction mapping is soft continuous.
**Proof.** Let $\tilde{x}_\lambda \in SP(X)$ be any soft point and an arbitrary $\tilde{\varepsilon} \tilde{>} \bar{0}$ is given. If we choose $d(\tilde{x}_\lambda, \tilde{y}_\mu) \tilde{<} \tilde{\delta} \tilde{<} \tilde{\varepsilon}$, then since $(f,\varphi)$ is soft contraction mapping we have
$$\tilde{d}((f,\varphi)(\tilde{x}_\lambda),(f,\varphi)(\tilde{y}_\mu)) \tilde{\leq} \alpha.\tilde{d}(\tilde{x}_\lambda,\tilde{y}_\mu) < \alpha.\delta \tilde{\leq} \tilde{\varepsilon}$$
and so $(f,\varphi)$ is soft continuous.

**Theorem 4.8.** Let $(\tilde{X}, \tilde{d}, E)$ be a soft complete metric space. If the mapping $(f,\varphi):(\tilde{X},\tilde{d},E) \to (\tilde{X},\tilde{d},E)$ is a soft contraction mapping on a complete soft metric space, then there exists a unique soft point $\tilde{x}_\lambda \in SP(X)$ such that $(f,\varphi)(\tilde{x}_\lambda) = \tilde{x}_\lambda$.
**Proof.** Let $\tilde{x}_\lambda^0$ be any soft point in $SP(X)$. Set
$$\tilde{x}_{\lambda_1}^1 = (f,\varphi)(\tilde{x}_\lambda^0) = \left(f(\tilde{x}_\lambda^0)\right)_{\varphi(\lambda)}, \tilde{x}_{\lambda_2}^2 = ((f,\varphi)(\tilde{x}_{\lambda_1}^1)) = \left(f^2(\tilde{x}_\lambda^0)\right)_{\varphi^2(\lambda)},\ldots$$
$$,\tilde{x}_{\lambda_{n+1}}^{n+1} = ((f,\varphi)(\tilde{x}_{\lambda_n}^n)) = \left(f^{n+1}(\tilde{x}_\lambda^0)\right)_{\varphi^{n+1}(\lambda)},\ldots$$

We have
$$\tilde{d}(\tilde{x}_{\lambda_{n+1}}^{n+1},\tilde{x}_{\lambda_n}^n) = \tilde{d}\left((f,\varphi)(\tilde{x}_{\lambda_n}^n),(f,\varphi)(\tilde{x}_{\lambda_{n-1}}^{n-1})\right) \tilde{\leq} \alpha.\tilde{d}(\tilde{x}_{\lambda_n}^n,\tilde{x}_{\lambda_{n-1}}^{n-1}) \tilde{\leq} \alpha^2.\tilde{d}(\tilde{x}_{\lambda_{n-1}}^{n-1},\tilde{x}_{\lambda_{n-2}}^{n-2}) \tilde{\leq} \ldots \tilde{\leq} \alpha^n.\tilde{d}(\tilde{x}_{\lambda_1}^1,\tilde{x}_\lambda^0).$$
So for $n > m$
$$\tilde{d}(\tilde{x}_{\lambda_n}^n,\tilde{x}_{\lambda_m}^m) \tilde{\leq} \tilde{d}(\tilde{x}_{\lambda_n}^n,\tilde{x}_{\lambda_{n-1}}^{n-1}) + \tilde{d}(\tilde{x}_{\lambda_{n-1}}^{n-1},\tilde{x}_{\lambda_{n-2}}^{n-2}) + \cdots + \tilde{d}(\tilde{x}_{\lambda_{m+1}}^{m+1},\tilde{x}_{\lambda_m}^m)$$
$$\tilde{\leq} (\alpha^{n-1} + \alpha^{n-2} + \cdots + \alpha^m)\tilde{d}(\tilde{x}_{\lambda_1}^1,\tilde{x}_\lambda^0)$$
$$\tilde{\leq} \frac{\alpha^m}{1-\alpha}\tilde{d}(\tilde{x}_{\lambda_1}^1,\tilde{x}_\lambda^0).$$

We get $\tilde{d}(\tilde{x}_{\lambda_n}^n,\tilde{x}_{\lambda_m}^m) \tilde{\leq} \frac{\alpha^m}{1-\alpha}\tilde{d}(\tilde{x}_{\lambda_1}^1,\tilde{x}_\lambda^0)$. This implies $\tilde{d}(\tilde{x}_{\lambda_n}^n,\tilde{x}_{\lambda_m}^m) \to \bar{0}$ $(n, m \to \infty)$. Hence $\{\tilde{x}_{\lambda_n}^n\}$ is a soft Cauchy sequence, by the completeness of $\tilde{X}$, there is $\tilde{x}_\lambda^* \tilde{\in} \tilde{X}$ such that $\tilde{x}_{\lambda_n}^n \to \tilde{x}_\lambda^*$ $(n \to \infty)$. Since

$$\tilde{d}\big((f,\varphi)(\tilde{x}_\lambda^*),\tilde{x}_\lambda^*\big) \tilde{\leq} \tilde{d}\big((f,\varphi)(\tilde{x}_{\lambda_n}^n),(f,\varphi)(\tilde{x}_\lambda^*)\big)+\tilde{d}\big((f,\varphi)(\tilde{x}_{\lambda_n}^n),\tilde{x}_\lambda^*\big)$$

$$\tilde{\leq} \alpha\tilde{d}\big(\tilde{x}_{\lambda_n}^n,\tilde{x}_\lambda^*\big)+\tilde{d}\big(\tilde{x}_{\lambda_{n+1}}^{n+1},\tilde{x}_\lambda^*\big),$$

$$\tilde{d}\big((f,\varphi)(\tilde{x}_\lambda^*),\tilde{x}_\lambda^*\big) \tilde{\leq} \alpha\Big(\alpha\tilde{d}\big(\tilde{x}_{\lambda_n}^n,\tilde{x}_\lambda^*\big)+\tilde{d}\big(\tilde{x}_{\lambda_{n+1}}^{n+1},\tilde{x}_\lambda^*\big)\Big)\to \bar{0}.$$

Hence $\tilde{d}\big((f,\varphi)(\tilde{x}_\lambda^*),\tilde{x}_\lambda^*\big)\to \bar{0}$. This implies $(f,\varphi)(\tilde{x}_\lambda^*)=\tilde{x}_\lambda^*$. So the soft point $\tilde{x}_\lambda^*$ is a fixed soft point of the mapping $(f,\varphi)$.

Now if $\tilde{y}_\mu^*$ is another fixed soft point of $(f,\varphi)$, then

$$\tilde{d}(\tilde{x}_\lambda^*,\tilde{y}_\mu^*)=\tilde{d}\big((f,\varphi)(\tilde{x}_\lambda^*),(f,\varphi)(\tilde{y}_\mu^*)\big)\tilde{\leq} \alpha\tilde{d}(\tilde{x}_\lambda^*,\tilde{y}_\mu^*).$$

Hence for $\tilde{\alpha}<\bar{1}$ $\tilde{d}(\tilde{x}_\lambda^*,\tilde{y}_\mu^*)=\bar{0}\Rightarrow \tilde{x}_\lambda^*=\tilde{y}_\mu^*$. Therefore the fixed soft point of $(f,\varphi)$ is unique.

**Theorem 4.9.** Let $(\tilde{X},\tilde{d},E)$ be a soft complete metric space. Suppose the soft mapping $(f,\varphi):(\tilde{X},\tilde{d},E)\to(\tilde{X},\tilde{d},E)$ satisfies the soft contractive condition

$$\tilde{d}\big((f,\varphi)(\tilde{x}_\lambda),(f,\varphi)(\tilde{y}_\mu)\big)\tilde{\leq}\tilde{\alpha}\Big[\tilde{d}\big((f,\varphi)(\tilde{x}_\lambda),\tilde{x}_\lambda\big)+\tilde{d}\big((f,\varphi)(\tilde{y}_\mu),\tilde{y}_\mu\big)\Big],\text{ for all } \tilde{x}_\lambda,\tilde{y}_\mu\tilde{\in}\tilde{X},$$

where $\tilde{\alpha}\in\left[\bar{0},\dfrac{\bar{1}}{2}\right)$ is a soft constant. Then $(f,\varphi)$ has a unique fixed soft point in $\tilde{X}$.

**Proof.** Choose $\tilde{x}_\lambda^0$ be any soft point in $SP(X)$. Set

$$\tilde{x}_{\lambda_1}^1=(f,\varphi)(\tilde{x}_\lambda^0)=\big(f(\tilde{x}_\lambda^0)\big)_{\varphi(\lambda)},\ \tilde{x}_{\lambda_2}^2=\big((f,\varphi)(\tilde{x}_{\lambda_1}^1)\big)=\big(f^2(\tilde{x}_\lambda^0)\big)_{\varphi^2(\lambda)},\dots$$

$$,\tilde{x}_{\lambda_{n+1}}^{n+1}=\big((f,\varphi)(\tilde{x}_{\lambda_n}^n)\big)=\big(f^{n+1}(\tilde{x}_\lambda^0)\big)_{\varphi^{n+1}(\lambda)},\dots$$

We have

$$\tilde{d}(\tilde{x}_{\lambda_{n+1}}^{n+1},\tilde{x}_{\lambda_n}^n)=\tilde{d}\big((f,\varphi)(\tilde{x}_{\lambda_n}^n),(f,\varphi)(\tilde{x}_{\lambda_{n-1}}^{n-1})\big)\tilde{\leq}\tilde{\alpha}\Big[\tilde{d}\big((f,\varphi)(\tilde{x}_{\lambda_n}^n),\tilde{x}_{\lambda_n}^n\big)+\tilde{d}\big((f,\varphi)(\tilde{x}_{\lambda_{n-1}}^{n-1}),\tilde{x}_{\lambda_{n-1}}^{n-1}\big)\Big]$$

So

$$=\tilde{\alpha}\Big[\tilde{d}\big(\tilde{x}_{\lambda_{n+1}}^{n+1},\tilde{x}_{\lambda_n}^n\big)+\tilde{d}\big(\tilde{x}_{\lambda_n}^n,\tilde{x}_{\lambda_{n-1}}^{n-1}\big)\Big].$$

$$\tilde{d}(\tilde{x}_{\lambda_{n+1}}^{n+1},\tilde{x}_{\lambda_n}^n)\tilde{\leq}\dfrac{\tilde{\alpha}}{\bar{1}-\tilde{\alpha}}\tilde{d}(\tilde{x}_{\lambda_n}^n,\tilde{x}_{\lambda_{n-1}}^{n-1})=h\tilde{d}(\tilde{x}_{\lambda_n}^n,\tilde{x}_{\lambda_{n-1}}^{n-1}),$$

where $\tilde{h}=\dfrac{\tilde{\alpha}}{\bar{1}-\tilde{\alpha}}$. For $n>m$,

$$\tilde{d}(\tilde{x}_{\lambda_n}^n,\tilde{x}_{\lambda_m}^m)\tilde{\leq}\tilde{d}(\tilde{x}_{\lambda_n}^n,\tilde{x}_{\lambda_{n-1}}^{n-1})+\tilde{d}(\tilde{x}_{\lambda_{n-1}}^{n-1},\tilde{x}_{\lambda_{n-2}}^{n-2})+\cdots+\tilde{d}(\tilde{x}_{\lambda_{m+1}}^{m+1},\tilde{x}_{\lambda_m}^m)$$

$$\tilde{\leq}(\tilde{h}^{n-1}+\tilde{h}^{n-2}+\cdots+\tilde{h}^m)\tilde{d}(\tilde{x}_{\lambda_1}^1,\tilde{x}_\lambda^0)$$

$$\tilde{\leq}\dfrac{\tilde{h}^m}{\bar{1}-\tilde{h}}\tilde{d}(\tilde{x}_{\lambda_1}^1,\tilde{x}_\lambda^0).$$

We get $\tilde{d}(\tilde{x}_{\lambda_n}^n,\tilde{x}_{\lambda_m}^m)\tilde{\leq}\dfrac{\tilde{h}^m}{\bar{1}-\tilde{h}}\tilde{\alpha}\tilde{d}(\tilde{x}_{\lambda_1}^1,\tilde{x}_\lambda^0)$. This implies $\tilde{d}(\tilde{x}_{\lambda_n}^n,\tilde{x}_{\lambda_m}^m)\to\bar{0}$ $(n,m\to\infty)$. Hence $\{\tilde{x}_{\lambda_n}^n\}$ is a soft Cauchy sequence. By the completeness of $\tilde{X}$, there is $\tilde{x}_\lambda^*\tilde{\in}\tilde{X}$ such that $\tilde{x}_{\lambda_n}^n\to\tilde{x}_\lambda^*$ $(n\to\infty)$. Since

$$\tilde{d}\big((f,\varphi)(\tilde{x}_\lambda^*),\tilde{x}_\lambda^*\big) \tilde{\leq} \tilde{d}\big((f,\varphi)(\tilde{x}_{\lambda_n}^n),(f,\varphi)(\tilde{x}_\lambda^*)\big) + \tilde{d}\big((f,\varphi)(\tilde{x}_{\lambda_n}^n),\tilde{x}_\lambda^*\big)$$

$$\tilde{\leq} \tilde{\alpha}\Big[\tilde{d}\big((f,\varphi)(\tilde{x}_{\lambda_n}^n),\tilde{x}_{\lambda_n}^n\big) + \tilde{d}\big((f,\varphi)(\tilde{x}_{\lambda_n}^n),\tilde{x}_\lambda^*\big)\Big] + \tilde{d}\big(\tilde{x}_{\lambda_{n+1}}^{n+1},\tilde{x}_\lambda^*\big)$$

$$\tilde{\leq} \frac{\bar{1}}{\bar{1}-\tilde{\alpha}}\Big[\tilde{\alpha}\tilde{d}\big((f,\varphi)(\tilde{x}_{\lambda_n}^n),\tilde{x}_{\lambda_n}^n\big) + \tilde{d}\big(\tilde{x}_{\lambda_{n+1}}^{n+1},\tilde{x}_\lambda^*\big)\Big],$$

$$\tilde{d}\big((f,\varphi)(\tilde{x}_\lambda^*),\tilde{x}_\lambda^*\big) \tilde{\leq} \tilde{\alpha}\frac{\bar{1}}{\bar{1}-\tilde{\alpha}}\Big(\tilde{\alpha}\tilde{d}\big(\tilde{x}_{\lambda_{n+1}}^{n+1},\tilde{x}_{\lambda_n}^n\big) + \tilde{d}\big(\tilde{x}_{\lambda_{n+1}}^{n+1},\tilde{x}_\lambda^*\big)\Big) \to \bar{0}.$$

Hence $\tilde{d}\big((f,\varphi)(\tilde{x}_\lambda^*),\tilde{x}_\lambda^*\big) \to \bar{0}$. This implies $(f,\varphi)(\tilde{x}_\lambda^*) = \tilde{x}_\lambda^*$. So the soft point $\tilde{x}_\lambda^*$ is a fixed soft point of the mapping $(f,\varphi)$.

Now if $\tilde{y}_\mu^*$ is another fixed soft point of $(f,\varphi)$, then

$$\tilde{d}(\tilde{x}_\lambda^*,\tilde{y}_\mu^*) = \tilde{d}\big((f,\varphi)(\tilde{x}_\lambda^*),(f,\varphi)(\tilde{y}_\mu^*)\big) \tilde{\leq} \tilde{\alpha}\Big[\tilde{d}\big((f,\varphi)(\tilde{x}_\lambda^*),\tilde{x}_\lambda^*\big) + \tilde{d}\big((f,\varphi)(\tilde{y}_\mu^*),\tilde{y}_\mu^*\big)\Big] = \bar{0}.$$

Hence for $\tilde{\alpha}<\bar{1}$ $\tilde{d}(\tilde{x}_\lambda^*,\tilde{y}_\mu^*) = \bar{0} \Rightarrow \tilde{x}_\lambda^* = \tilde{y}_\mu^*$. Therefore the fixed soft point of $(f,\varphi)$ is unique.

**Theorem 4.10.** Let $(\tilde{X},\tilde{d},E)$ be a soft complete metric space. Suppose the soft mapping $(f,\varphi):(\tilde{X},\tilde{d},E) \to (\tilde{X},\tilde{d},E)$ satisfies the contractive condition

$$\tilde{d}\big((f,\varphi)(\tilde{x}_\lambda),(f,\varphi)(\tilde{y}_\mu)\big) \tilde{\leq} \tilde{\alpha}\Big[\tilde{d}\big((f,\varphi)(\tilde{x}_\lambda),\tilde{y}_\mu\big) + \tilde{d}\big((f,\varphi)(\tilde{y}_\mu),\tilde{x}_\lambda\big)\Big], \text{ for all } \tilde{x}_\lambda,\tilde{y}_\mu \tilde{\in} \tilde{X}$$

where $\tilde{\alpha} \in \left[\bar{0},\frac{\bar{1}}{2}\right)$ is a constant. Then $(f,\varphi)$ has a unique fixed soft point in $\tilde{X}$.

**Proof.** Choose $\tilde{x}_\lambda^0$ be any soft point in $SP(X)$. Set

$$\tilde{x}_{\lambda_1}^1 = (f,\varphi)(\tilde{x}_\lambda^0) = \big(f(\tilde{x}_\lambda^0)\big)_{\varphi(\lambda)}, \tilde{x}_{\lambda_2}^2 = \big((f,\varphi)(\tilde{x}_{\lambda_1}^1)\big) = \big(f^2(\tilde{x}_\lambda^0)\big)_{\varphi^2(\lambda)}, \ldots$$

$$,\tilde{x}_{\lambda_{n+1}}^{n+1} = \big((f,\varphi)(\tilde{x}_{\lambda_n}^n)\big) = \big(f^{n+1}(\tilde{x}_\lambda^0)\big)_{\varphi^{n+1}(\lambda)}, \ldots$$

We have

$$\tilde{d}(\tilde{x}_{\lambda_{n+1}}^{n+1},\tilde{x}_{\lambda_n}^n) = \tilde{d}\big((f,\varphi)(\tilde{x}_{\lambda_n}^n),(f,\varphi)(\tilde{x}_{\lambda_{n-1}}^{n-1})\big) \tilde{\leq} \tilde{\alpha}\Big[\tilde{d}\big((f,\varphi)(\tilde{x}_{\lambda_n}^n),\tilde{x}_{\lambda_{n-1}}^{n-1}\big) + \tilde{d}\big((f,\varphi)(\tilde{x}_{\lambda_{n-1}}^{n-1}),\tilde{x}_{\lambda_n}^n\big)\Big]$$

$$\tilde{\leq} \tilde{\alpha}\Big[\tilde{d}\big(\tilde{x}_{\lambda_{n+1}}^{n+1},\tilde{x}_{\lambda_n}^n\big) + \tilde{d}\big(\tilde{x}_{\lambda_n}^n,\tilde{x}_{\lambda_{n-1}}^{n-1}\big)\Big].$$

So,

$$\tilde{d}(\tilde{x}_{\lambda_{n+1}}^{n+1},\tilde{x}_{\lambda_n}^n) \tilde{\leq} \frac{\tilde{\alpha}}{\bar{1}-\tilde{\alpha}}\tilde{d}(\tilde{x}_{\lambda_n}^n,\tilde{x}_{\lambda_{n-1}}^{n-1}) = \tilde{h}\tilde{d}(\tilde{x}_{\lambda_n}^n,\tilde{x}_{\lambda_{n-1}}^{n-1}),$$

where $\tilde{h} = \frac{\tilde{\alpha}}{\bar{1}-\tilde{\alpha}}$. For $n>m$,

$$\tilde{d}(\tilde{x}_{\lambda_n}^n,\tilde{x}_{\lambda_m}^m) \tilde{\leq} \tilde{d}(\tilde{x}_{\lambda_n}^n,\tilde{x}_{\lambda_{n-1}}^{n-1}) + \tilde{d}(\tilde{x}_{\lambda_{n-1}}^{n-1},\tilde{x}_{\lambda_{n-2}}^{n-2}) + \cdots + \tilde{d}(\tilde{x}_{\lambda_{m+1}}^{m+1},\tilde{x}_{\lambda_m}^m)$$

$$\tilde{\leq} (\tilde{h}^{n-1} + \tilde{h}^{n-2} + \cdots + \tilde{h}^m)\tilde{d}(\tilde{x}_{\lambda_1}^1,\tilde{x}_\lambda^0)$$

$$\tilde{\leq} \frac{\tilde{h}^m}{\bar{1}-\tilde{h}}\tilde{d}(\tilde{x}_{\lambda_1}^1,\tilde{x}_\lambda^0).$$

We get $\tilde{d}(\tilde{x}^n_{\lambda_n}, \tilde{x}^m_{\lambda_m}) \tilde{\leq} \dfrac{\tilde{h}^m}{1-\tilde{h}} \tilde{\alpha} \tilde{d}(\tilde{x}^1_{\lambda_1}, \tilde{x}^0_\lambda)$. This implies $\tilde{d}(\tilde{x}^n_{\lambda_n}, \tilde{x}^m_{\lambda_m}) \to \overline{0}$ $(n,m \to \infty)$. Hence $\{\tilde{x}^n_{\lambda_n}\}$ is a soft Cauchy sequence. By the completeness of $\tilde{X}$, there is $\tilde{x}^*_\lambda \tilde{\in} \tilde{X}$ such that $\tilde{x}^n_{\lambda_n} \to \tilde{x}^*_\lambda$ $(n \to \infty)$. Since

$$\tilde{d}\big((f,\varphi)(\tilde{x}^*_\lambda), \tilde{x}^*_\lambda\big) \tilde{\leq} \tilde{d}\big((f,\varphi)(\tilde{x}^n_{\lambda_n}), (f,\varphi)(\tilde{x}^*_\lambda)\big) + \tilde{d}\big((f,\varphi)(\tilde{x}^n_{\lambda_n}), \tilde{x}^*_\lambda\big)$$

$$\tilde{\leq} \tilde{\alpha}\Big[\tilde{d}\big((f,\varphi)(\tilde{x}^*_\lambda), \tilde{x}^n_{\lambda_n}\big) + \tilde{d}\big((f,\varphi)(\tilde{x}^n_{\lambda_n}), \tilde{x}^*_\lambda\big)\Big] + \tilde{d}\big(\tilde{x}^{n+1}_{\lambda_{n+1}}, \tilde{x}^*_\lambda\big)$$

$$\tilde{\leq} \tilde{\alpha}\Big[\tilde{d}\big((f,\varphi)(\tilde{x}^*_\lambda), \tilde{x}^*_\lambda\big) + \tilde{d}\big(\tilde{x}^n_{\lambda_n}, \tilde{x}^*_\lambda\big) + \tilde{d}\big(\tilde{x}^{n+1}_{\lambda_{n+1}}, \tilde{x}^*_\lambda\big)\Big] + \tilde{d}\big(\tilde{x}^{n+1}_{\lambda_{n+1}}, \tilde{x}^*_\lambda\big)$$

$$\tilde{\leq} \dfrac{\overline{1}}{\overline{1}-\tilde{\alpha}}\Big[\tilde{\alpha}\big(\tilde{d}\big(\tilde{x}^n_{\lambda_n}, \tilde{x}^*_\lambda\big) + \tilde{d}\big(\tilde{x}^{n+1}_{\lambda_{n+1}}, \tilde{x}^*_\lambda\big)\big) + \tilde{d}\big(\tilde{x}^{n+1}_{\lambda_{n+1}}, \tilde{x}^*_\lambda\big)\Big],$$

$$\tilde{d}\big((f,\varphi)(\tilde{x}^*_\lambda), \tilde{x}^*_\lambda\big) \tilde{\leq} \tilde{\alpha} \dfrac{\tilde{\alpha}}{\overline{1}-\tilde{\alpha}}\Big[\tilde{\alpha}\big(\tilde{d}\big(\tilde{x}^n_{\lambda_n}, \tilde{x}^*_\lambda\big) + \tilde{d}\big(\tilde{x}^{n+1}_{\lambda_{n+1}}, \tilde{x}^*_\lambda\big)\big) + \tilde{d}\big(\tilde{x}^{n+1}_{\lambda_{n+1}}, \tilde{x}^*_\lambda\big)\Big] \to \overline{0}.$$

Hence $\tilde{d}\big((f,\varphi)(\tilde{x}^*_\lambda), \tilde{x}^*_\lambda\big) \to \overline{0}$. This implies $(f,\varphi)(\tilde{x}^*_\lambda) = \tilde{x}^*_\lambda$. So the soft point $\tilde{x}^*_\lambda$ is a fixed soft point of the mapping $(f,\varphi)$.

Now if $\tilde{y}^*_\mu$ is another fixed soft point of $(f,\varphi)$, then

$$\tilde{d}\big(\tilde{x}^*_\lambda, \tilde{y}^*_\mu\big) = \tilde{d}\big((f,\varphi)(\tilde{x}^*_\lambda), (f,\varphi)(\tilde{y}^*_\mu)\big) \tilde{\leq} \tilde{\alpha}\Big[\tilde{d}\big((f,\varphi)(\tilde{x}^*_\lambda), \tilde{y}^*_\mu\big) + \tilde{d}\big((f,\varphi)(\tilde{y}^*_\mu), \tilde{x}^*_\lambda\big)\Big] = \overline{0}.$$

Hence for $\tilde{\alpha} \tilde{<} \overline{1}$ $\tilde{d}\big(\tilde{x}^*_\lambda, \tilde{y}^*_\mu\big) = \overline{0} \Rightarrow \tilde{x}^*_\lambda = \tilde{y}^*_\mu$. Therefore the fixed soft point of $(f,\varphi)$ is unique.

**Proposition 4.11.** Let $(\tilde{X}, \tilde{d}, E)$ be a soft metric space. If $(f,\varphi):(\tilde{X},\tilde{d},E) \to (\tilde{X},\tilde{d},E)$ is a soft contraction mapping, then the mapping $f_\lambda:(X,d_\lambda) \to (X, d_{\varphi(\lambda)})$ is contraction mapping, for all $\lambda \in E$.

The following example shows that converse of the Proposition 4.11 does not hold.

**Example 4.12.** Let $E = \mathbb{R}$ be parameter set and $X = \mathbb{R}^2$. Consider usual metrics on this sets and define soft metric on $\tilde{X}$ by $\tilde{d}(\tilde{x}_\lambda, \tilde{y}_\mu) = |\lambda - \mu| + d(x,y)$. Then if we define the soft mapping $(f,\varphi):(\tilde{X},\tilde{d},E) \to (\tilde{X},\tilde{d},E)$ as follows

$$(f,\varphi)(x_\lambda) = \left(\dfrac{1}{2}x\right)_{3\lambda},$$

then

$$\tilde{d}((f,\varphi)(0,1)_2, (f,\varphi)(1,0)_1) = \tilde{d}\left((0,\dfrac{1}{2})_6, (\dfrac{1}{2},0)_3\right) = 3 + \dfrac{\sqrt{2}}{2}$$

$$\tilde{d}\big((0,1)_2, (1,0)_1\big) = 1 + \sqrt{2}.$$

Since $3 + \dfrac{\sqrt{2}}{2} > 1 + \sqrt{2}$, we see that the soft mapping $(f,\varphi)$ is not a soft contraction mapping. But the mapping $f:(X, d_\lambda) \to (X, d_{3\lambda})$ is a contraction mapping, for all $\lambda \in E$.

**Corallary 4.13.** Let $E$ be a parameter set and $X$ be a set. By using the given metrics defined on these sets, we can form a soft metric. If $(f, \varphi)$ is a soft contraction mapping on the soft space then $f$ or $\varphi$ may not be a contraction mapping.

The following example justifies the Corallary 4.13.

**Example 4.14.** Let $(\tilde{\mathbb{R}}, \tilde{d}, E)$ be a soft metric space with the following metrics

$$d(x, y) = |x - y|, \quad d_1(x, y) = \min\{|x - y|, 1\} \text{ and}$$

$$\tilde{d}(x_\lambda, y_\mu) = \frac{1}{2} d_1(\lambda, \mu) + d(x, y)$$

where $E = [1, \infty)$. Let the functions $\varphi : [1, \infty) \to [1, \infty)$ and $f : \mathbb{R} \to \mathbb{R}$ are defined as follows
$\varphi(x) = x + \frac{1}{x}$ and $f(x) = \frac{1}{5} x$ respectively.

Here, it is obvious that the conditions of contraction mapping are hold for the composite function

$$(f, \varphi) : (\tilde{\mathbb{R}}, \tilde{d}, E) \to (\tilde{\mathbb{R}}, \tilde{d}, E) .$$

We want to show that $(f, \varphi)$ is a soft contraction mapping, whereas the function $\varphi(x) = x + \frac{1}{x}$ is not a contraction mapping with the defined metric $d_1(x, y) = \min\{|x - y|, 1\}$.

$$\tilde{d}((f, \varphi)(\tilde{x}_\lambda), (f, \varphi)(\tilde{y}_\mu)) = \tilde{d}\left(\left(\frac{1}{5} x\right)_{\lambda + \frac{1}{\lambda}}, \left(\frac{1}{5} y\right)_{\mu + \frac{1}{\mu}}\right)$$

$$= \frac{1}{2} d_1\left(\lambda + \frac{1}{\lambda}, \mu + \frac{1}{\mu}\right) + \frac{1}{5} |x - y|$$

$$= \frac{1}{2} \min\left\{\left|\lambda + \frac{1}{\lambda} - \mu + \frac{1}{\mu}\right|, 1\right\} + \frac{1}{5} |x - y|$$

$$= \frac{1}{2} \min\left\{|\lambda - \mu|\left|1 - \frac{1}{\lambda \mu}\right|, 1\right\} + \frac{1}{5} |x - y|$$

$$\leq \frac{1}{2} \min\{|\lambda - \mu|, 1\} + \frac{1}{5} |x - y|$$

$$\leq \frac{1}{2} d_1(\lambda, \mu) + \frac{1}{5} d(x, y)$$

$$\leq \frac{3}{4} (d_1(\lambda, \mu) + d(x, y)),$$

which shows that $(f, \varphi)$ is a soft contraction mapping.

## 5. Conclusion.

In the present work, we have continued to investigate the properties of soft metric spaces. We also introduce soft continuous mappings. Later we prove some fixed point theorems of soft contractive mappings on soft metric spaces. We hope that the findings in this paper will help researcher enhance and promote the further study on soft metric spaces to carry out a general framework for their applications in real life.

## References

[1] M.I. Ali, F. Feng, X. Liu, W.K. Min and M. Shabir, On some new operations in soft set theory, Comput. Math. Appl.49(2005) 1547-1553.

[2] D.Chen, The parameterization reduction of soft sets and its applications, Comput. Math. Appl. 49 (2005) 757-763.

[3] Sujoy Das and S. K. Samanta, Soft real sets, soft real numbers and their properties, J. Fuzzy Math. 20 (3) (2012) 551-576.

[4] Sujoy Das and S. K. Samanta, Soft metric, Annals of Fuzzy Mathematics and Informatics, 6(1) (2013) 77-94.

[5] Sujoy Das and S. K. Samanta, On soft metric spaces, J. Fuzzy Math, accepted.

[6] C. Gunduz (Aras), A. Sonmez, H. Çakallı, On Soft Mappings, arXiv:1305.4545v1 [math.GM], 16 May 2013.

[7] S. Hussain and B. Ahmad, Some properties of soft topological spaces, Computers and Math. with Applications, 62(2011), 4058-4067.

[8] P. K. Maji, A.R.Roy, R. Biswas, An application of soft sets in a decision making problem, Comput. Math. Appl.44 (2002) 1077-1083.

[9] P.K.Maji, R.Biswas, A.R.Roy, Soft set theory, Comput. Math. Appl.45 (2003) 555-562.

[10] P. Majumdar and S. K. Samanta, On soft mappings, Comput. Math. Appl. 60 (2010) 2666-2672.

[11] D. Molodtsov, Soft set-theory-first results, Comput. Math. Appl.37(1999) 19-31.

[12] M. Shabir and M. Naz, On soft topological spaces, Comput. Math. Appl.61(2011) 1786-1799.

[13] S. Bayramov, C. Gunduz(Aras), Soft locally compact and soft paracompact spaces, Journal of Mathematics and System Science (accepted)

[14] B.E. Rhoades, A comparison of various definition of contractive mappings, Trans. Amer. Math. Soc. 266 (1977) 257-290

[15] H. Long-Guang, Z. Xian, Cone metric spaces and fixed point theorems of contractive mappings, J. Math. Anal. Appl. 332 (2007) 1468-1476